\documentclass[12pt]{amsart}
\usepackage[left=1.5in,right=1.5in,top=1.5in,bottom=1.5in]{geometry}
\usepackage{amsmath}
\usepackage{amssymb}
\usepackage{fancyhdr}
\usepackage{hyperref}
\usepackage{graphicx}

\hypersetup{
  colorlinks,
  filecolor = blue,
  citecolor = blue,
  urlcolor  = blue,
  linkcolor = blue
}

\newcommand{\refeq}[1] {
  \hyperref[#1]{(\ref*{#1})}
}
\numberwithin{equation}{section}

\title[Spherical Triangle]
{Whitney Forms for Spherical Triangles I: \\
The Euler, Cagnoli, and Tuynman \\
Area Formulas,
Barycentric Coordinates, \\
and Construction
with the Exterior Calculus}

\author{David W. Fillmore}
\address{Tech-X Corporation, Boulder, Colorado}
\email{fillmore@txcorp.com}
\author{Jay P. Fillmore}
\address{Department of Mathematics,
University of California, San Diego (retired)}

\begin{document}
\maketitle
\begin{center}
\today
\end{center}

\begin{abstract}
We establish the equivalence
of the Tuynman midpoint area formula
for a spherical triangle
to the classical area formulas of Euler and of Cagnoli.
The derivation also yields a variant
of the Cagnoli formula in terms of the medial triangle.

We introduce the three barycentric coordinates of a
point within the spherical triangle
as area ratios which sum to unity.
The barycentric coordinates are the Whitney 0-forms,
scalar functions over the domain of the triangle
associated with each vertex.
We then construct,
by exterior differentiation of the barycentric coordinates,
succinct expressions for the Whitney 1-forms
associated with each geodesic side, or great circle arc.

The Euler formula, in conjunction with that of Tuynman,
facilitates the differentiation of a triangular area
with respect to the position of a vertex.
As both the Euler and Tuynman formulas
may be expressed naturally in terms of the position vectors
of the vertices on an embedded sphere,
the Whitney constructions may be done
in terms of vector-valued forms
and without recourse to a particular projection or coordinate chart.

Finally, we exhibit the Whitney 2-form of the triangle,
which must be the product of a scalar function
and the area 2-form of the sphere.
We find an expression for this scalar function
in terms of determinants of $3 \times 3$
matrices built from the vertex position vectors.
It is a rational function in the Cartesian coordinates of a point.
By construction it must be invariant under cyclic permutations
of the three vertices,
though it is not manifestly so.
Also by construction the Whitney 2-form
must integrate to one over the triangle.

We speculate that these exactly integrable rational functions
associated with spherical triangles may be relevant
to the development of numerical quadrature schemes on the sphere,
as well as be of inherent interest in their own right,
and that the spherical Whitney forms may have potential application
in the formulation of discrete exterior calculus
and finite element spaces intrinsic to the sphere.
\end{abstract}

\newpage
\section{Area Formulas for a Spherical Triangle}
Consider an oriented spherical triangle $\triangle ABC$
on the unit sphere $\mathbb{S}^2$ embedded in $\mathbb{R}^3$.
Let $A$, $B$ and $C$ be the unit vectors in $\mathbb{R}^3$
from the center of $\mathbb{S}^2$ to the vertices of $\triangle ABC$.
In the common nomenclature of spherical trigonometry,
$a$, $b$ and $c$ denote the lengths of the sides
opposite the vertices $A$, $B$ and $C$, respectively.
These spherical lengths are angles between
the unit vectors in $\mathbb{R}^3$,
so that $\cos a = B \cdot C$,
$\cos b = C \cdot A$, and $\cos c = A \cdot B$.
Let $M(A, B, C)$ be the $3 \times 3$ matrix formed
from the components of $A$, $B$, and $C$ as columns,
and let $S$ be the oriented (signed) area of $\triangle ABC$.
The sine of half the oriented area is
given by a formula recently
established by Tuynman \cite{Tuynman2013}:
\begin{equation}
\label{equation:sphere-triangle-area-tuynman}
\sin \frac{S}{2}
= \frac{\det M(A, B, C)}
{\sqrt{2 \left( 1 + A \cdot B \right)
\left( 1 + B \cdot C \right)
\left( 1 + C \cdot A \right)
}} \ .
\end{equation}
Note that the factors in the denominator
may be expressed in terms of the side lengths through
$1 + A \cdot B = 1 + \cos c = 2 \cos^2 \frac{c}{2}$,
and similarly for the other factors,
so that equation (\ref{equation:sphere-triangle-area-tuynman})
may be written as
\begin{equation}
\label{equation:sphere-triangle-area-intermediate}
\sin \frac{S}{2}
= \frac{\sqrt{\det \left( M^{\mathrm{T}} M \right)}}
{4 \cos \frac{a}{2} \cos \frac{b}{2} \cos \frac{c}{2}} \ .
\end{equation}
In equation (\ref{equation:sphere-triangle-area-intermediate})
$M(A, B, C)$ has been abbreviated simply as $M$
and we have replaced $\det M$
with $\sqrt{\left( \det M \right)^2}
= \sqrt{\det \left( M^{\mathrm{T}} M \right)}$
in anticipation of the algebraic manipulations to follow.
In equation (\ref{equation:sphere-triangle-area-tuynman})
and in many to follow,
none of the dot products between the vertex unit vectors
can be allowed to take the value $-1$,
that is, for the vertices to be antipodal.
Thus we must exclude triangles with antipodal vertices
from further consideration,
or otherwise resolve the indeterminant form $\frac{0}{0}$
with l'H\^{o}pital's rule.
For simplicity of exposition, we choose to do the former.
Moreover, we will take only the positive root of $\det \left( M^{\mathrm{T}} M \right)$.
This is tantamount to a requirement that all triangles under consideration
have a positive orientation,
so that {\it all} determinants that will appear in this note are strictly
greater than zero.
To emphasize: $S > 0$ and $\det M(A, B, C) > 0$ throughout.
The Tuynman formula is closely related to
the classical area formulas of Euler and of Cagnoli
(see, for instance, the book by Donnay \cite{Donnay1945}).
To establish the equivalence
of these formulas we form the
symmetric matrix product
\begin{equation}
\label{equation:matrix-product}
M^{\mathrm{T}} M
= \left( \begin{array}{ccc}
A \cdot A & A \cdot B & A \cdot C \\
B \cdot A & B \cdot B & B \cdot C \\
C \cdot A & C \cdot B & C \cdot C
\end{array} \right)
= \left( \begin{array}{ccc}
1 & \cos c & \cos b \\
\cos c & 1 & \cos a \\
\cos b & \cos a & 1
\end{array} \right)
\end{equation}
and find its determinant to be
\begin{equation}
\label{equation:matrix-product-determinant}
\det \left( M^{\mathrm{T}} M \right)
= 1 - \cos^2 a - \cos^2 b - \cos^2 c
+ 2 \cos a \cos b \cos c \ .
\end{equation}
The Euler formula relates
the {\it cosine} of half the area to the side lengths
while the Cagnoli formula,
as does equation (\ref{equation:sphere-triangle-area-intermediate}),
relates the {\it sine} of half the area to the lengths.
From the trigonometric identity
$\cos \frac{S}{2} = \sqrt{1 - \sin^2 \frac{S}{2}}$
we have
\begin{equation}
\label{equation:intermediate-step}
\cos \frac{S}{2}
= \frac{ \sqrt{
\left(
4 \cos \frac{a}{2} \cos \frac{b}{2} \cos \frac{c}{2}
\right)^2
- \det \left( M^{\mathrm{T}} M \right)}}
{4 \cos \frac{a}{2} \cos \frac{b}{2} \cos \frac{c}{2}} \ .
\end{equation}
The radicand
in equation (\ref{equation:intermediate-step}) is
\begin{equation}
\begin{aligned}
& \ 2 \left( 1 + \cos a \right)
\left( 1 + \cos b \right)
\left( 1 + \cos c \right)
- \det \left( M^{\mathrm{T}} M \right) \\
&=
2 \left( 1 + \cos a + \cos b + \cos c \right)
+ 2 \left( \cos a \cos b
+ \cos b \cos c + \cos c \cos a \right) \\
&+ 2 \cos a \cos b \cos c
- 1 + \cos^2 a + \cos^2 b + \cos^2 c
- 2 \cos a \cos b \cos c \\
&= \left( 1 + \cos a + \cos b + \cos c \right)^2
\end{aligned}
\end{equation}
from which the Euler formula follows:
\begin{equation}
\label{equation:sphere-triangle-area-euler}
\cos \frac{S}{2}
= \frac{1 + \cos a + \cos b + \cos c}
{4 \cos \frac{a}{2} \cos \frac{b}{2} \cos \frac{c}{2}} \ .
\end{equation}
Equation (\ref{equation:sphere-triangle-area-euler})
may be recast with vector dot products in $\mathbb{R}^3$,
similar to equation (\ref{equation:sphere-triangle-area-tuynman}),
as
\begin{equation}
\label{equation:sphere-triangle-area-euler-vector}
\cos \frac{S}{2}
= \frac{1 + A \cdot B + B \cdot C + C \cdot A}
{\sqrt{2 \left( 1 + A \cdot B \right)
\left( 1 + B \cdot C \right)
\left( 1 + C \cdot A \right)
}} \ .
\end{equation}
In a similar fashion one may
begin with equations
(\ref{equation:sphere-triangle-area-intermediate})
and (\ref{equation:matrix-product-determinant}),
and through the judicious use of standard trigonometric identities,
arrive at the classical formula of Cagnoli:
\begin{equation}
\label{equation:sphere-triangle-area-cagnoli}
\sin \frac{S}{2}
= \frac{\sqrt{\sin s \sin (s - a) \sin (s - b) \sin (s - c)}}
{2 \cos \frac{a}{2} \cos \frac{b}{2} \cos \frac{c}{2}} \ ,
\end{equation}
where $s = \frac{a + b + c}{2}$ is the semi-perimeter
of $\triangle ABC$.
As an excursion, let us rediscover a variant
of the Cagnoli formula in terms of the medial triangle
$\triangle DEF$ of triangle $\triangle ABC$.
As discussed by Tuynman \cite{Tuynman2013},
equation (\ref{equation:sphere-triangle-area-tuynman})
takes the simple form
\begin{equation}
\label{equation:sphere-triangle-area-tuynman-midpoint}
\sin \frac{S}{2} = \det M(D, E, F)
\end{equation}
when expressed in terms
of the side midpoint unit vectors
$D = \frac{A + B}{\sqrt{2(1 + A \cdot B)}}$,
$E = \frac{B + C}{\sqrt{2(1 + B \cdot C)}}$
and
$F = \frac{C + A}{\sqrt{2(1 + C \cdot A)}}$.
The midpoint formula
(\ref{equation:sphere-triangle-area-tuynman-midpoint})
follows immediately from the properties
of $\det M(D, E, F)$:
linearity in each column
and anti-symmetry with respect to the interchange of columns.
From equation (\ref{equation:sphere-triangle-area-tuynman-midpoint})
it is clear that $\left| \sin \frac{S}{2} \right| \le 1$,
as required.
Let $d$, $e$ and $f$ denote the lengths of the sides
opposite the vertices in $\triangle DEF$.
By direct comparison
of equation (\ref{equation:sphere-triangle-area-intermediate})
versus (\ref{equation:sphere-triangle-area-cagnoli})
we have, without further effort,
the Cagnoli equivalent of
(\ref{equation:sphere-triangle-area-tuynman-midpoint}):
\begin{equation}
\label{equation:sphere-triangle-area-cagnoli-midpoint}
\sin \frac{S}{2}
= 2 \sqrt{\sin t \sin (t - d) \sin (t - e) \sin (t - f)},
\end{equation}
where $t = \frac{d + e + f}{2}$ is the semi-perimeter
of the medial triangle $\triangle DEF$.
For completeness we also present an equivalent
of the Euler formula
(\ref{equation:sphere-triangle-area-euler}):
\begin{equation}
\label{equation:sphere-triangle-area-euler-midpoint}
\cos \frac{S}{2}
= \sqrt{\cos^2 d + \cos^2 e + \cos^2 f - 2 \cos d \cos e \cos f} \ .
\end{equation}
Let us conclude this section with a few remarks.
In the subsequent development
we will not make any use of the interior angles,
conventionally denoted $\alpha$, $\beta$ and $\gamma$,
associated with the vertices $A$, $B$ and $C$, respectively.
These angles are also the dihedral angles between
the planes spanned, in the case of $\gamma$ for instance,
the basis sets \{ $A$, $C$ \} and \{ $B$, $C$ \}.
The dihedral product (a notion valid in any dimension)
between these planes is
$ \sin a \sin b \cos \gamma
= (A \cdot B) (C \cdot C) - (A \cdot C) (B \cdot C)
= \cos c - \cos a \cos b $,
the law of cosines for a spherical triangle.
The area of a spherical triangle is known to be the same as its
angular excess $S = \alpha + \beta + \gamma - \pi$.
While, through the law of cosines,
each interior angle can be related to the arc-cosine
of some combination of vector dot products,
we have yet to explore this avenue.
\section{Differentiation of the Area of a Spherical Triangle}
Let us take a point $X$ interior to or on the boundary of
$\triangle ABC$
and construct the edges $XA$,
$XB$ and $XC$.
This results in a partition of $\triangle ABC$
into three {\it positively} oriented triangles
$\triangle ABX$, $\triangle BCX$ and $\triangle CAX$
with a common vertex $X$.
Denote their respective area functions by
$S_C(X)$, $S_A(X)$ and $S_B(X)$.
Should $X$ coincide with any of the vertices
$A$, $B$ or $C$ we must recover the undivided triangle so that
$S_A(A) = S$, $S_B(B) = S$ and $S_C(C) = S$.
The three area functions must sum to the total area of $\triangle ABC$
\begin{equation}
S_A(X) + S_B(X) + S_C(X) = S \ .
\end{equation}
Define $M_A(X) = M(B, C, X)$ and
$f_A : \mathbb{S}^2 \rightarrow \mathbb{R}$ as
\begin{equation}
\label{equation:f}
f_A\left( X \right)
= 2 \left( 1 + B \cdot C \right)
\left( 1 + B \cdot X \right)
\left( 1 + C \cdot X \right)
\end{equation}
so that equations
(\ref{equation:sphere-triangle-area-euler-vector})
and
(\ref{equation:sphere-triangle-area-tuynman})
may be written concisely as
\begin{equation}
\label{equation:sphere-triangle-area-euler-brief}
\cos \frac{S_A(X)}{2} = \frac{1 + B \cdot C + (B + C) \cdot X}{\sqrt{f_A(X)}}
\end{equation}
and
\begin{equation}
\label{equation:sphere-triangle-area-tuynman-brief}
\sin \frac{S_A(X)}{2} = \frac{\det M_A(X)}{\sqrt{f_A(X)}} \ .
\end{equation}
Similar expressions for $M_B(X), f_B(X), S_B(X)$
and $M_C(X), f_C(X), S_C(X)$
are obtained through cyclic permutation of the vertices $A, B, C$.
Again, in the formation of the $3 \times 3$ matrices
we must keep the columns ordered as $M(A, B, X)$, $M(B, C, X)$, and $M(C, A, X)$.
We wish to find the exterior derivative $d$ of $S_A(X)$
in terms of the vector-valued 1-form $dX$.
Upon application of $d$ to equation
(\ref{equation:sphere-triangle-area-euler-brief})
we have
\begin{equation}
- \frac{1}{2} \sin \frac{S_A}{2} dS_A
= \frac{(B + C) \cdot dX}{\sqrt{f_A}}
- \frac{1}{2} \cos \frac{S_A}{2} \frac{df_A}{f_A}
\end{equation}
or, from equation (\ref{equation:sphere-triangle-area-tuynman-brief}),
\begin{equation}
- \frac{1}{2} \det M_A \ dS_A
= (B + C) \cdot dX
- \frac{1}{2} \left[ 1 + B \cdot C + \left( B + C \right) \cdot X \right]
\frac{df_A}{f_A} \ .
\end{equation}
From equation (\ref{equation:f}) we have
\begin{equation}
\frac{df_A}{f_A}
= \frac{B \cdot dX}{1 + B \cdot X}
+ \frac{C \cdot dX}{1 + C \cdot X} \ .
\end{equation}
After some rearrangement our final result for $dS_A$ is
\begin{equation}
- \det M_A \ dS_A
= \left[ 1 - \frac{C \cdot (X + B)}{1 + B \cdot X} \right] B \cdot dX
+ \left[ 1 - \frac{B \cdot (X + C)}{1 + C \cdot X} \right] C \cdot dX \ ,
\end{equation}
with similar formulas for $dS_B$ and $dS_C$
obtained through cyclic permutation of $A$, $B$ and $C$.
\section{The Whitney Forms of a Spherical Triangle}
Let $\lambda_A(X) = \frac{S_A(X)}{S}$,
$\lambda_B(X) = \frac{S_B(X)}{S}$
and $\lambda_C(X) = \frac{S_C(X)}{S}$
be the barycentric coordinates of the point $X$
interior to or on the boundary of $\triangle ABC$,
such that $\lambda_A(X) + \lambda_B(X) + \lambda_C(X) = 1$.
Should $X$ reside on the side $AB$, $\lambda_C(X) = 0$
and $\lambda_A(X) + \lambda_B(X) = 1$.
Should $X$ coincide with the vertex $A$,
$\lambda_B(X) = 0$, $\lambda_C(X) = 0$
and $\lambda_A(X) = 1$.
Similar statements hold for the other two sides
and the other two vertices.
The $\lambda_i(X)$ are the three Whitney 0-forms
associated with the vertices $i = A, B, C$ of the triangle.
Recall that the three Whitney 1-forms
associated with the sides of the triangle
are defined through exterior differentiation as
\begin{equation}
\lambda_{ij} = \lambda_i \ d\lambda_j - \lambda_j \ d\lambda_i
\end{equation}
for $ij$ = $AB$, $BC$ and $CA$ \cite{Whitney1957}.
%
%
Let $\gamma_{ij}$ be the great circle arc
that joins the vertices $i$ and $j$.
The Whitney forms are defined such that
$\int_{\gamma_{AB}} \lambda_{AB} = 1$,
$\int_{\gamma_{BC}} \lambda_{AB} = 0$,
$\int_{\gamma_{CA}} \lambda_{AB} = 0$,
and so forth.
Let us introduce the vector function
$F_{BC} : \mathbb{S}^2 \rightarrow T_{\mathbb{S}^2} \mathbb{R}^3$ as
\begin{equation}
\label{equation:F_definition}
F_{BC}(X) = \left[ 1 - \frac{C \cdot (X + B)}{1 + B \cdot X} \right] B
+ \left[ 1 - \frac{B \cdot (X + C)}{1 + C \cdot X} \right] C
\end{equation}
so that
\begin{equation}
- \det M_A \ dS_A = F_{BC} \cdot dX
\end{equation}
or
\begin{equation}
\label{equation:d_lambda_A}
d\lambda_A = - \frac{F_{BC}}{S \ \det M_A} \cdot dX \ ,
\end{equation}
and in a similar manner
\begin{equation}
\label{equation:d_lambda_B}
d\lambda_B = - \frac{F_{CA}}{S \ \det M_B} \cdot dX
\end{equation}
and
\begin{equation}
\label{equation:d_lambda_C}
d\lambda_C = - \frac{F_{AB}}{S \ \det M_C} \cdot dX \ .
\end{equation}
Note that equation (\ref{equation:d_lambda_A}) for $d \lambda_A$
is indeterminant of order $\frac{0}{0}$
in the limits $X \rightarrow B$ and $X \rightarrow C$,
and likewise for $d \lambda_B$ and $d \lambda_C$.
That this is so can be seen
from the definition (\ref{equation:F_definition})
and that fact that $X \cdot dX = 0$.
Since the denominators are all linear in $X$,
all of these limits must be finite.
To elaborate somewhat on our formalism, $T \mathbb{R}^3$
is the tangent bundle to $\mathbb{R}^3$
(the manifold of tangent spaces at each point)
and $T_{\mathbb{S}^2} \mathbb{R}^3$
is the tangent bundle $T \mathbb{R}^3$
restricted to the surface of the unit sphere.
We assert that the $F_{ij}(X)$ should be in $T_X \mathbb{R}^3$,
as their genesis was through differentiation at the point $X$.
This is just a fancy way to say that the
vector $F_{ij}(X)$ should be thought of as anchored at $X$,
and that a vector field in $T_{\mathbb{S}^2} \mathbb{R}^3$
should be thought of as a field of vectors all anchored
to the surface of $\mathbb{S}^2$ but which point in any direction
in $\mathbb{R}^3$.
These careful distinctions are not necessary for actual calculations.
To continue with the Whitney forms,
the 1-form $\lambda_{AB}$, for instance, is
\begin{equation}
\label{equation:whitney-1-form}
\lambda_{AB}
= \frac{1}{S} \left( \frac{\lambda_B}{\det M_A} F_{BC}
- \frac{\lambda_A}{\det M_B} F_{AC} \right) \cdot dX \ ,
\end{equation}
where we have also used the symmetry $F_{CA} = F_{AC}$
to keep the formula explicitly anti-symmetric in $A$ and $B$.
As the functions $\lambda_i(X)$ must
contain either an arc-cosine or an arc-sine,
depending on the choice of area formula,
we have not found any further simplifications thus far.
We turn now to the Whitney 2-form
\begin{equation}
\label{equation:d_lambda_ABC}
\lambda_{ABC} = d \lambda_{AB}
= 2 \ d\lambda_A \wedge d\lambda_B \ .
\end{equation}
Since $\lambda_A + \lambda_B + \lambda_C = 1$,
one may in fact use any side,
as $d\lambda_A \wedge d\lambda_B
= d\lambda_B \wedge d\lambda_C = d\lambda_C \wedge d\lambda_A$.
From equations
(\ref{equation:d_lambda_A}) and (\ref{equation:d_lambda_B})
the Whitney 2-form is
\begin{equation}
\label{equation:whitney-2-form}
\lambda_{ABC}
= \frac{2}{S^2 \det \left( M_A M_B \right)}
\left( F_{BC} \cdot dX \right)
\wedge \left( F_{CA} \cdot dX \right) \ .
\end{equation}
The Whitney 2-form
(\ref{equation:d_lambda_ABC})
has a very special property:
$\label{equation:whitney-2-form-integral}
\int_{\triangle ABC} \lambda_{ABC} = 1$ \cite{Whitney1957}.
This is easy to see when the integral is done in barycentric coordinates:
\begin{equation}
2 \int_{\triangle ABC} d\lambda_A \wedge d\lambda_B
= 2 \int_0^1 d\lambda_A
\int_0^{1 - \lambda_A} d\lambda_B
= 1 \ .
\end{equation}
Note that the lower limit of integration for $\lambda_B$
is $0$, when $\lambda_C = 1 - \lambda_A$,
and that the upper limit is $1 - \lambda_A$, when $\lambda_C = 0$.
$S$ is the area of $\triangle ABC$
and so might better be written as an explicit
function of the vertices as $S(A, B, C)$.
Now any 2-form on $\mathbb{S}^2$ is the product of some
function $g(X) : \mathbb{S}^2 \rightarrow \mathbb{R}$
and the area 2-form $\mathrm{vol}^2$, so that
\begin{equation}
\lambda_{ABC}(X) = \frac{\omega_{ABC}(X)}{S(A, B, C)} \mathrm{vol}^2 \ .
\end{equation}
We wish to find $\omega_{ABC}(X)$.
To do this we jump back into $\mathbb{R}^3$
and introduce a radial coordinate $r$.
A position vector $x$ in $\mathbb{R}^3$
is just $x = r X$
where $r$ is its distance from
the center of $\mathbb{S}^2$.
Now the volume 3-form in $\mathbb{R}^3$
should be $\mathrm{vol}^3 = r^2 dr \wedge \mathrm{vol}^2$.
Since $x \cdot x = r^2$ we have $r \ dr = x \cdot dx$
or $dr = X \cdot dx$.
Also note that
$r \ dr \wedge \left( V \cdot dX \right)
= dr \wedge \left( V \cdot dx \right)$
for any vector $V$, so that
\begin{equation}
\begin{aligned}
r^2 dr \wedge \left( \frac{S \lambda_{ABC}}{\omega_{ABC}} \right)
&= \frac{2 r^2}{S \omega_{ABC} \det \left( M_A M_B \right)}
dr \wedge \left( F_{BC} \cdot dX \right)
\wedge \left( F_{CA} \cdot dX \right) \\
&= \frac{2}{S \omega_{ABC} \det \left( M_A M_B \right)}
\left( X \cdot dx \right)
\wedge \left( F_{BC} \cdot dx \right)
\wedge \left( F_{CA} \cdot dx \right) \ .
\end{aligned}
\end{equation}
After a moment of contemplation
one should realize that
\begin{equation}
\left( U \cdot dx \right)
\wedge \left( V \cdot dx \right)
\wedge \left( W \cdot dx \right)
= \det M(U, V, W) \ \mathrm{vol}^3
\end{equation}
for any three vectors $U, V, W$
in $\mathbb{R}^3$, ergo:
\begin{equation}
r^2 dr \wedge \left( \frac{S \lambda_{ABC}}{\omega_{ABC}} \right)
= \frac{2 \det M \left(F_{BC}, F_{CA}, X \right) \mathrm{vol}^3}
{S \omega_{ABC}
\det \left[ M \left( B, C, X \right)
M \left( C, A, X \right) \right]} \ .
\end{equation}
Now since $ r^2 dr \wedge \left( \frac{S \lambda_{ABC}}{\omega_{ABC}} \right)
= \mathrm{vol}^3$,
it must be that
\begin{equation}
\label{equation:rational_function}
\omega_{ABC}(X)
= \frac{2 \det M \left[F_{BC}(X), F_{CA}(X), X \right]}
{S \left( A, B, C \right) \det \left[ M \left( B, C, X \right)
M \left( C, A, X \right) \right]} \ .
\end{equation}
From equation (\ref{equation:F_definition})
one may observe
that in general $\omega_{ABC}(X)$
appears to be the sum of three rational functions
in the Cartesian coordinates of $X$,
each of which is the ratio of a third order polynomial
to a fourth order polynomial.
By the Whitney construction $\omega_{ABC}$
must be invariant under cyclic permutation of $A, B, C$,
however a direct proof from equation
(\ref{equation:rational_function}) has so far eluded us.
One might also notice that $\omega_{ABC}$ inherits
indeterminant limits from the Whitney 1-forms,
which can be seen all to be finite.
Plots of $\omega_{ABC}$ for selected triangles
can be found on pages 10 through 12.
In article II of this series,
we will have a much closer look
at the visualization of $\omega_{ABC}$,
as well as the Whitney 1-forms,
for various triangles.
We will also work out in more detail
special examples,
such as the $90^\circ - 90^\circ - 90^\circ$
and the $\alpha - 90^\circ - 90^\circ$
isosceles triangles.
To conclude our development,
we delve a bit further into the notion of the vector-valued 1-form $dX$.
Let $e_i$ for $i = 1, 2$ be an orthonormal
vector frame on the tangent bundle $T\mathbb{S}^2$,
the {\it rep\`ere mobile} of Cartan
\cite{Cartan1937, Cartan1945}
and $\omega^i$ its dual 1-form basis
in the cotangent bundle $T^*\mathbb{S}^2$
such that $\omega^i(e_j) = \delta^i_{\ j}$.
In the formalism of Cartan
$dX = e_i \omega^i$,
with summation implied over repeated lower and upper indices.
This states that the vector-valued 1-form
$dX$ is a map
$dX : T_{\mathbb{S}^2} \mathbb{R}^3 \rightarrow T\mathbb{S}^2$.
That is to say,
for any vector field $V$ in $T_{\mathbb{S}^2} \mathbb{R}^3$
the interior product
$V \lrcorner \ dX = dX(V) = e_i V^i$
returns that part of $V$ tangential to $\mathbb{S}^2$.
This we already know,
as $X \cdot X = 1$ implies $X \cdot dX = 0$,
or to belabor the point,
$V \lrcorner \ \left( X \cdot dX \right)
= X \cdot \left( V \lrcorner \ dX \right) = 0$.
To make the ideas of Cartan more tangible,
let us work out explicitly the simple case
in which the standard spherical coordinates
$(\theta, \phi)$ are used to locate
a point $X$ on the unit sphere.
In this case
$X = \left( \begin{array}{ccc} \sin \theta \cos \phi &
\sin \theta \sin \phi & \cos \theta \end{array} \right)^{\mathrm{T}}$
and
\begin{equation}
\begin{aligned}
dX
&= \left( \begin{array}{c}
\cos \theta \cos \phi \\
\cos \theta \sin \phi \\
- \sin \theta
\end{array} \right) d\theta
+ \left( \begin{array}{c}
- \sin \theta \sin \phi \\
\sin \theta \cos \phi \\
0
\end{array} \right) d\phi \\
&= \left( \begin{array}{c}
\cos \theta \cos \phi \\
\cos \theta \sin \phi \\
- \sin \theta
\end{array} \right) \omega^{\theta}
+ \left( \begin{array}{c}
- \sin \phi \\
\cos \phi \\
0
\end{array} \right) \omega^{\phi} \ ,
\end{aligned}
\end{equation}
as $\omega^{\theta} = d\theta$
and $\omega^\phi = \sin \theta \ d\phi$
are an orthonormal basis on $T^*\mathbb{S}^2$
such that the metric tensor is
$g = \omega^{\theta} \otimes \omega^{\theta}
+ \omega^{\phi} \otimes \omega^{\phi}$.
But
\begin{equation}
\begin{aligned}
e_{\theta}
&= \cos \theta \cos \phi \ e_x
+ \cos \theta \sin \phi \ e_y
- \sin \theta \ e_z \ , \\
e_{\phi}
&= - \sin \phi \ e_x + \cos \phi \ e_y
\end{aligned}
\end{equation}
comprise the familiar orthonormal frame on $T\mathbb{S}^2$
for which $g(e_i, e_j) = \delta_{ij}$,
so that indeed
$dX = e_{\theta} \omega^{\theta} + e_{\phi} \omega^{\phi}$.
%

%
%
\newpage
\section{Conclusions and Future Directions}
In this paper we have found compact formulas
for the Whitney forms of a spherical triangle.
Given the plethora of identities in spherical trigonometry,
there are doubtless many alternative ways
to effect this construction.
Of all the alternatives that we have explored,
however, the method presented here is the most succinct.
There are numerous possible continuations
and applications of these developments.
First and foremost,
we wish to obtain a deeper understanding
of the rational functions on $\mathbb{R}^3$ that
relate a Whitney 2-form to the area 2-form.
The $\omega_{ABC}$ rationals have three poles on $\mathbb{S}^2$
that coincide with the antipodal points of the vertices,
and we would like to find a relationship between the extrema
and the vertices.
Perhaps one can approach this question through
the pullback of the Whitney forms,
via an inverse stereographic projection,
onto the extended complex plane $\mathbb{C} \cup \{ \infty \}$.
%
%

%
In the realm of partial differential equations,
Whitney forms are relevant
to the finite element and discrete exterior calculus
\cite{Desbrun2005}.
The spherical Whitney forms may have potential use
in the formulation of finite element spaces
intrinsic to the sphere.
In such a case, one expects to be confronted
with integrals such as
$\int_{\triangle} \lambda_{AB} \wedge \star \ \lambda_{BC}$
in various combinations
and for every triangle $\triangle$
in a geodesic triangulation of $\mathbb{S}^2$.
These integrals must be computed
through numerical quadrature,
and implementation
of algorithms for this purpose may be pursued
in future work.
Generalization to three-dimensional geometry
when spherical coordinate systems are employed
may also be of interest.
Here one could attempt to define
barycentric coordinates and find the Whitney forms
for a right triangular prism
in which the bottom and top faces
of the prism are spherical triangles
situated on spheres of different radii.
The Tuynman formula also holds
for triangles of the hyperbolic plane $\mathbb{H}^2$
embedded in three-dimensional
Minkowski space, that is, $\mathbb{R}^3$
endowed with a metric $\mathrm{diag}(-1, 1, 1)$.
The Whitney constructions certainly have
analogs on $\mathbb{H}^2$,
and we would like to see the $\omega_{ABC}$
rationals on the Poincar\'{e} disc.
Further questions naturally arise, such as:
does the Tuynman formula generalize to higher dimensions,
such as to the volume of a spherical tetrahedron in $\mathbb{S}^3$?
Does it generalize to two-dimensional de Sitter space?
The original development of the Tuynman formula
and its precursors occurred in the context
of stereographic projection and the extended complex plane
\cite{Tuynman2013}.
Other techniques in spherical trigonometry that exploit
the conformal nature of the stereographic projection
were developed by Ces\`aro
in the early part of last century \cite{Donnay1945}.
While the approach in this paper was strictly algebraic,
we intend to explore further the rich geometric aspects
of these problems in subsequent articles.
\newpage
\section*{}
\begin{figure}
\includegraphics[scale=0.4]{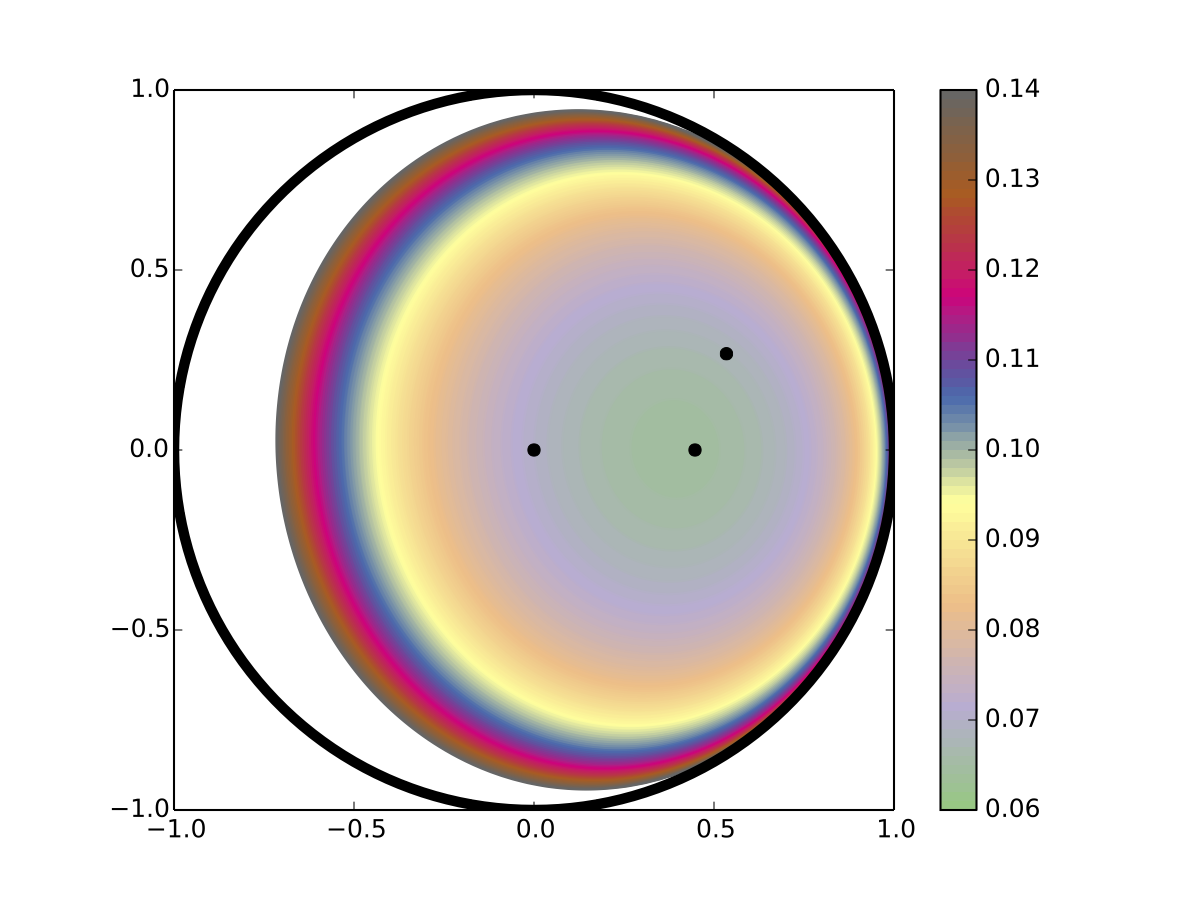}
\caption{$S \ \omega_{ABC}$ in the upper hemisphere $z > 0$ of $\mathbb{S}^2$,
projected by $(x, y, z) \mapsto (x, y)$.
The vertices
$A = \frac{e_x + 2 e_z}{\sqrt{2}},
B = \frac{2 e_x + e_y + 3 e_z}{\sqrt{14}},
C = e_z$
are indicated by block dots.}
\end{figure}
\begin{figure}
\includegraphics[scale=0.4]{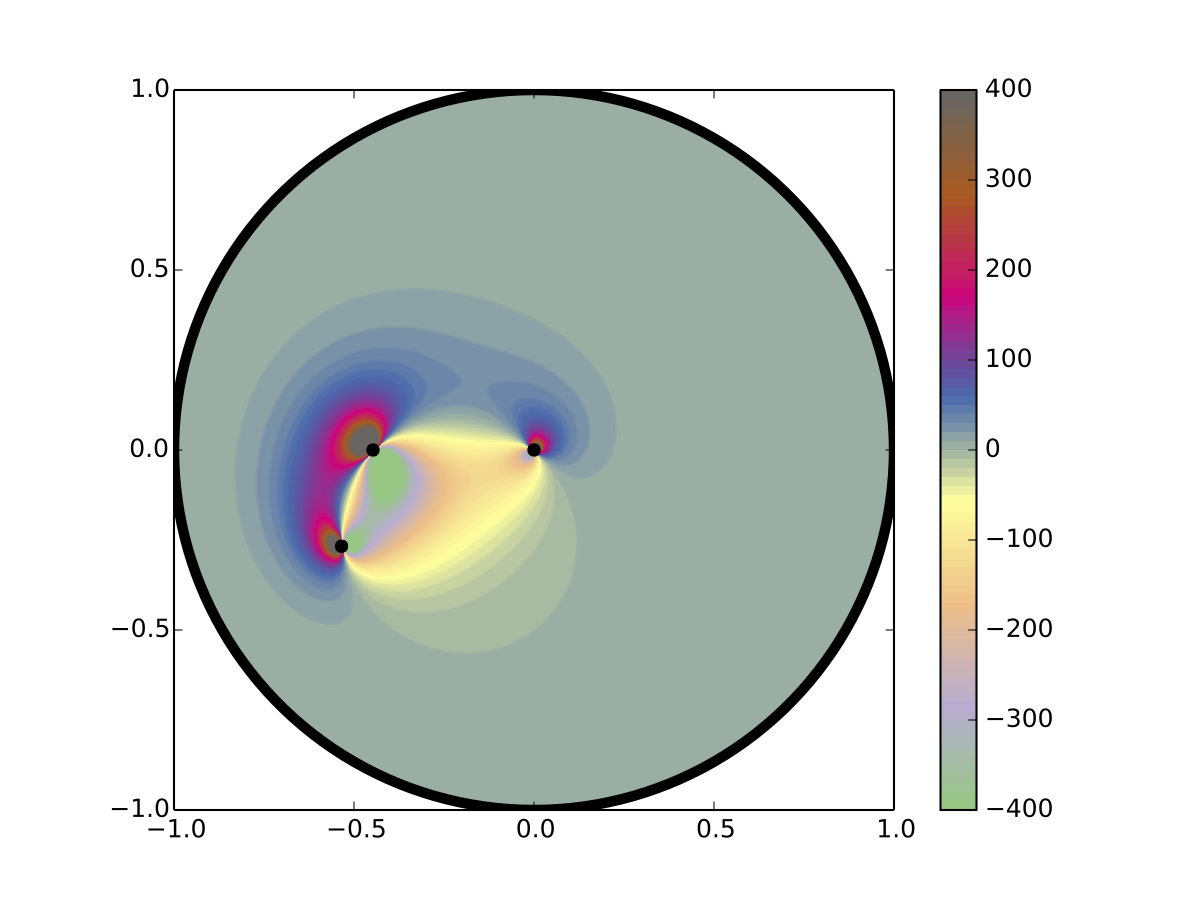}
\caption{$S \ \omega_{ABC}$ in the lower hemisphere $z < 0$ of $\mathbb{S}^2$,
projected by $(x, y, z) \mapsto (x, y)$.
The antipodal vertices $- A, - B, - C$ are indicated by block dots.}
\end{figure}
\clearpage

\begin{figure}
\includegraphics[scale=0.4]{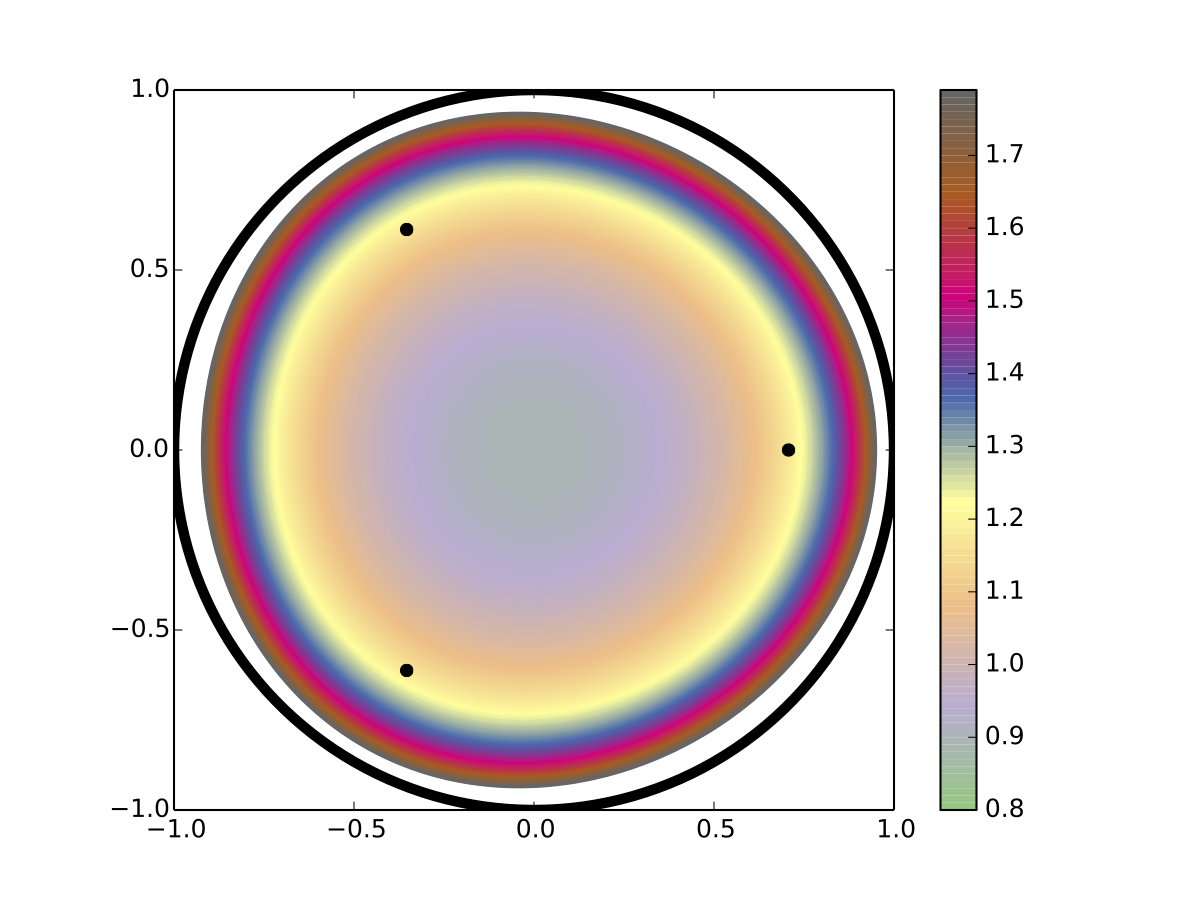}
\caption{$S \ \omega_{ABC}$ in the upper hemisphere $z > 0$ of $\mathbb{S}^2$,
projected by $(x, y, z) \mapsto (x, y)$.
The vertices of this equilateral triangle
are $A = \frac{e_x + e_z}{\sqrt{2}}$,
with $B$ and $C$ rotations of $A$
about the $z-$ axis by $\frac{2 \pi}{3}$ and $\frac{4 \pi}{3}$.
They are indicated by block dots.}
\end{figure}
\begin{figure}
\includegraphics[scale=0.4]{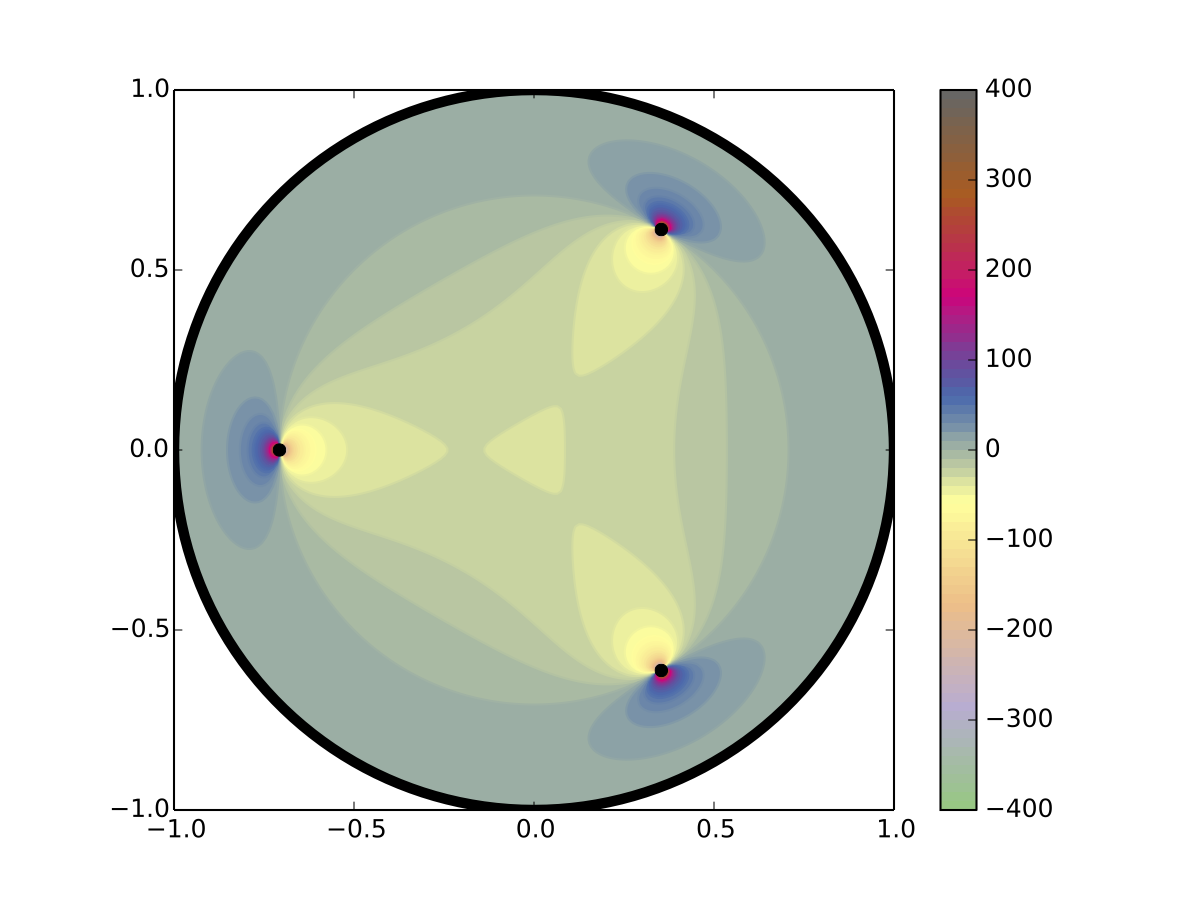}
\caption{$S \ \omega_{ABC}$ in the lower hemisphere $z < 0$ of $\mathbb{S}^2$,
projected by $(x, y, z) \mapsto (x, y)$.
The antipodal vertices $- A, - B, - C$ are indicated by block dots.}
\end{figure}
\clearpage

\section*{}
\begin{figure}
\includegraphics[scale=0.4]{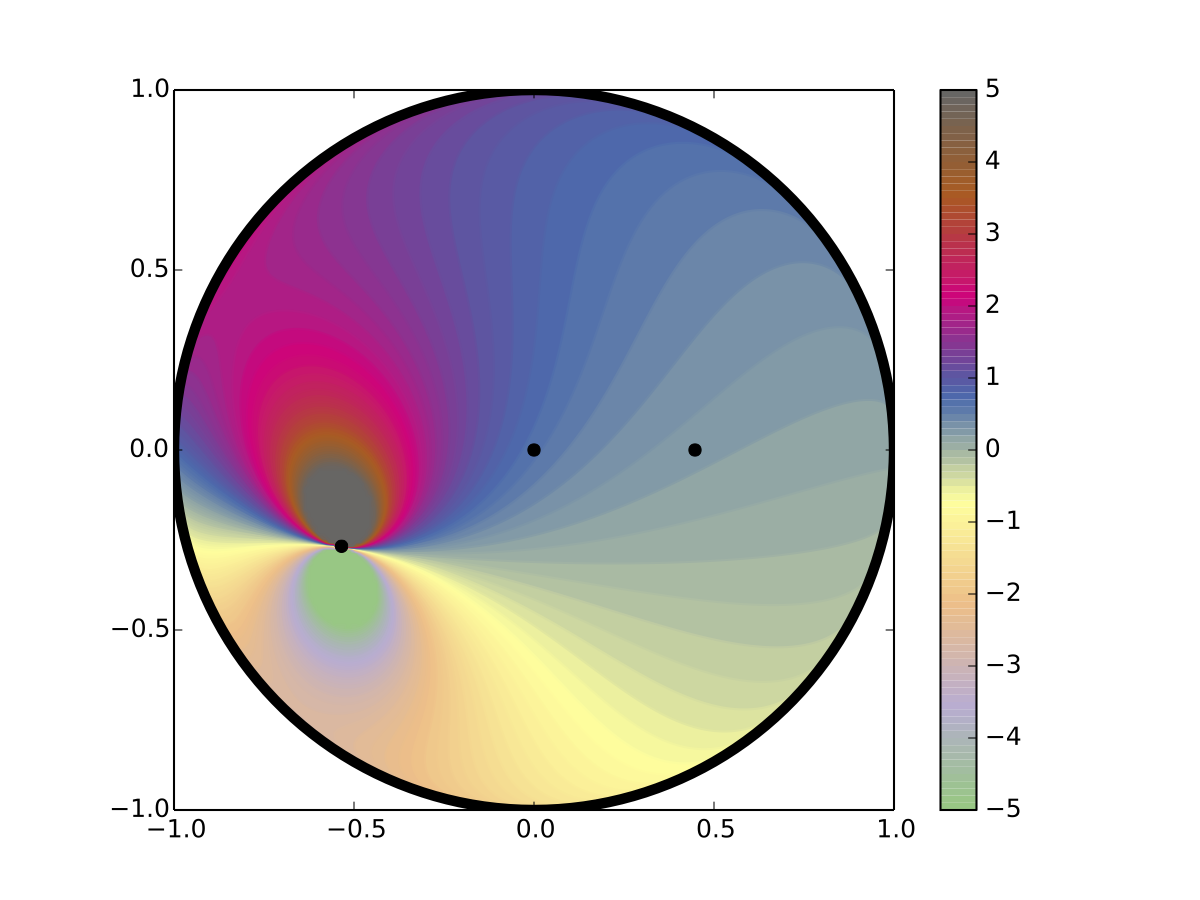}
\caption{$S \ \omega_{ABC}$ in the upper hemisphere $z > 0$ of $\mathbb{S}^2$,
projected by $(x, y, z) \mapsto (x, y)$.
The points
$A = \frac{e_x + 2 e_z}{\sqrt{2}},
- B = - \frac{2 e_x + e_y - 3 e_z}{\sqrt{14}},
C = e_z$
are indicated by block dots.}
\end{figure}
\begin{figure}
\includegraphics[scale=0.4]{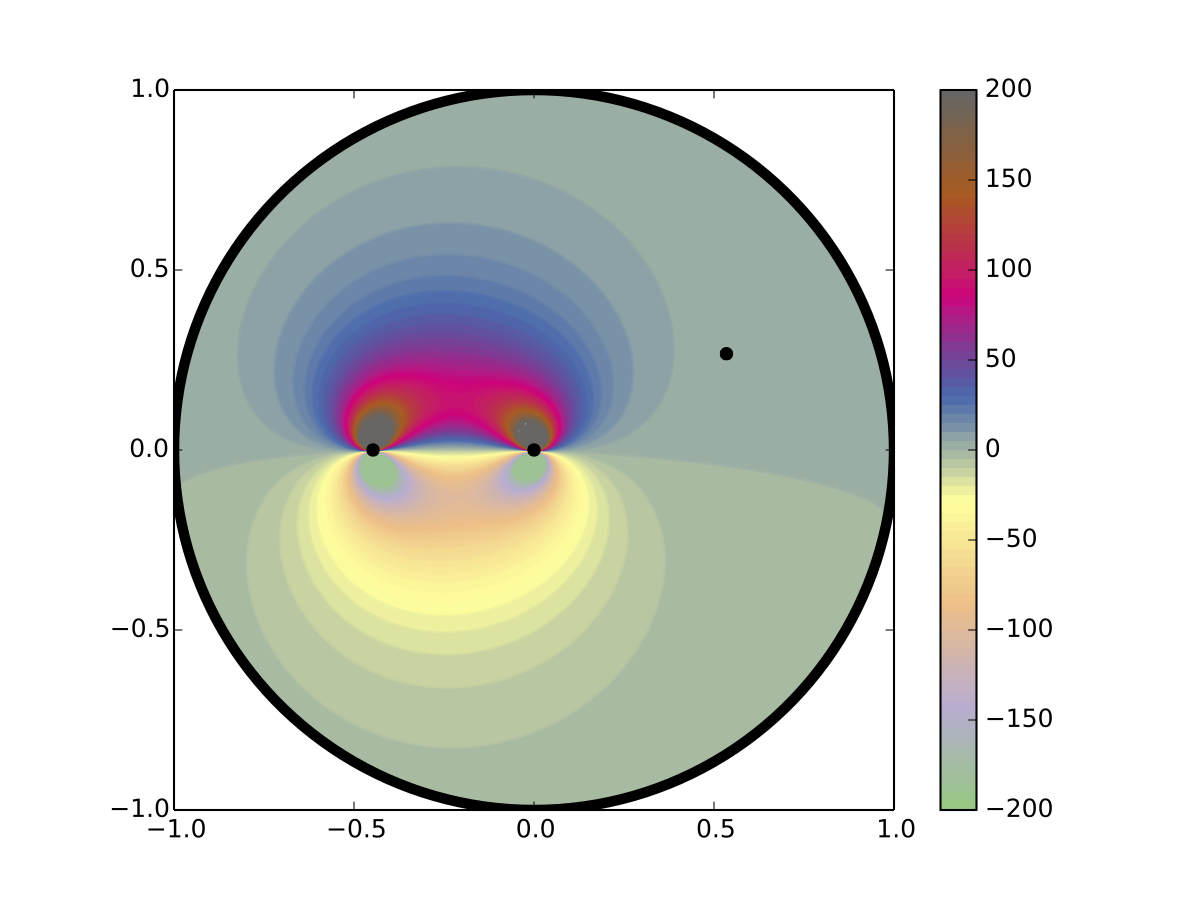}
\caption{$S \ \omega_{ABC}$ in the lower hemisphere $z < 0$ of $\mathbb{S}^2$,
projected by $(x, y, z) \mapsto (x, y)$.
The points $- A, B, - C$ are indicated by block dots.}
\end{figure}

\vspace{0.5in}
\section*{Acknowledgments}
We wish to thank Dr. John Thuburn,
Dr. Werner Bauer, and Dr. Colin Cotter
for helpful feedback and corrections.

The lead author Dr. David Fillmore
wishes to thank
Dr. Alexander Pletzer of Tech-X Corporation
for many illuminating discussions on
the applications of the exterior calculus,
and for partial support
under the DOE SBIR Phase II grant
\url{http://www.sbir.gov/sbirsearch/detail/390667},
{\it Visualizing Staggered Vector Fields}.
I also wish to thank my father and co-author,
Dr. Jay Fillmore,
who showed me the beauty in mathematics.
\end{document}